\input amstex
\documentstyle{amsppt}
\topmatter \magnification=\magstep1 \pagewidth{5.2 in}
\pageheight{6.7 in}
\abovedisplayskip=10pt \belowdisplayskip=10pt
\parskip=8pt
\parindent=5mm
\baselineskip=2pt
\title
  A note on $q$-Bernoulli numbers and polynomials
\endtitle
\author   Taekyun Kim    \endauthor

\affil{ {\it Institute of Science Education,\\
        Kongju National University, Kongju 314-701, S. Korea}}\endaffil
\thanks  2000 Mathematics Subject Classification:  11S80, 11B68, 11M99 .\endthanks
\abstract{ By using $q$-integration, we will give some integral
equation which are related to  the Barnes' multiple Bernoulli
numbers. The object of this paper is to give explicit
$q$-integral's formulae  which are related to Barnes' multiple
$q$-Bernoulli polynomials.
 }\endabstract
\rightheadtext{  $q$-Bernoulli numbers and Polynomials}
\leftheadtext{T. Kim}
\endtopmatter

\document

\define\Z{\Bbb Z_p}

\head \S 1. Introduction \endhead
 Let  $\Bbb Z,\,\Bbb Z_p,\,\Bbb Q_p $ and $\Bbb C_p$ will be  denoted by
  the ring of integers, the ring
of $p$-adic integers, the field of $p$-adic numbers and  the
completion of algebraic closure of $\Bbb Q_p$ and let $\nu_p$ be
the normalized exponential valuation of $\Bbb C_p$ with
$|p|_p=p^{-\nu_p(p)}=p^{-1}.$  When one talks of $q$-extension,
$q$ is variously considered as an indeterminate, a complex number
or a $p$-adic number.  If $ q \in \Bbb C_p ,$  then we  normally
assume $|q-1|_p < p^{-\frac1{p-1}},$ so that $q^x=\exp(x\log q)$
for $|x|_p \leq 1.$

In this paper, we use the notation:
$$[x]=[x:q]=\frac{1-q^x}{1-q}.$$
Hence,  $\lim_{q\rightarrow1}[x]=x$ for any $x$ with $|x|_p \leq 1
$ in the present $p$-adic case.

Let $d$ be a fixed integer and let $p$ be a fixed prime number. We
set
$$\align
&X=\varprojlim_N (\Bbb Z/dp^N\Bbb Z),\\
&X^*=\bigcup\Sb 0<a<dp\\ (a,p)=1\endSb a+dp\Bbb Z_p,\\
&a+dp^N\Bbb Z_p=\{x\in X\mid x\equiv a\pmod{dp^N}\},
\endalign$$
where $a\in \Bbb Z$ lies in $0\leq a<dp^N,$ cf. [2, 3, 4, 5, 6, 7,
8].

For any positive integer $N,$  we set $$\mu_q(a+dp^N
\Z)=\frac{q^a}{[dp^N]}$$  and this can be extended to a
distribution on $X.$ This distribution yields an integral for each
non-negative integer $m$:
$$ \int_{\Bbb Z_p} [x]^m d\mu_q(x)=\int_{X}[a]^{m}d\mu_q(a)=\frac1{(1-q)^m}
\sum_{i=0}^m\binom mi(-1)^i \frac{i+1}{[i+1]} .$$ The multiple
Barnes' Bernoulli polynomials were defined by
$$\left( \prod_{j=1}^r  \frac{w_j} {e^{w_j t} -1} \right) {t^r} e^{xt}=\sum_{n=0}^{\infty}
B_n^{(r)}(x|w_1 , w_2, \cdots , w_r )\frac{ t^n}{n!}, \tag 1 $$
for  each $w_j >0, 0<t<1, $ cf. [1].

The numbers $B_n^{(r)}(w_1 , w_2, \cdots , w_r)=B_n^{(r)}(0|w_1 ,
w_2, \cdots , w_r)$ are called the multiple Barnes' numbers.

Throughout this paper, we assume that $\alpha_1, \cdots , \alpha_k
$ are taken in the set of integers and let $ w \in \Bbb Z_p $.
Now, we can consider the multiple $q$-Bernoulli numbers by using
$q$-integral as follows:
$$ \beta_n^{(r)}(w, q |\alpha_1 , \alpha_2, \cdots ,
\alpha_r)=\undersetbrace\text{$r$
times}\to{\int_{\Z}\int_{\Z}\cdots\int_{\Z}} [w + \alpha_1x_1
+\cdots +\alpha_r x_r]^n d\mu_q(x_1)  \cdots d\mu_q(x_r),$$ and
$$ \beta_n^{(r)}( q |\alpha_1 , \alpha_2, \cdots ,
\alpha_r)=\undersetbrace\text{$r$
times}\to{\int_{\Z}\int_{\Z}\cdots\int_{\Z}} [ \alpha_1x_1 +\cdots
+\alpha_r x_r]^n d\mu_q(x_1)  \cdots d\mu_q(x_r).$$

In this paper, we construct the numbers $\beta_n^{(h,r)}(w, q
|\alpha_1 , \alpha_2, \cdots , \alpha_r) $ for $ h \in \Bbb Z $
which reduce  the multiple Barnes' Bernoulli numbers
$B_n^{(r)}(\alpha_1 , \alpha_2, \cdots , \alpha_r)$ as $q
\rightarrow 1$. Also, we give the new explicit formulas which are
related to these numbers.

 \head 2. An extension of Changhee $q$-Bernoulli numbers  \endhead
For $ h \in \Bbb Z $, we define the extension of Changhee
$q$-Bernoulli polynomials, numbers as follows:
$$
\aligned
 & \beta_n^{(h, r)}(w, q |\alpha_1 , \alpha_2, \cdots ,
 \alpha_r) \\
 & =
 \undersetbrace\text{$r$ times}\to{\int_{\Z} \cdots\int_{\Z}}
q^{\sum_{i=1}^{r} (h-i)x_i}[w + \alpha_1x_1 +\cdots +\alpha_r
x_r]^n d\mu_q(x_1) \cdots d\mu_q(x_r),
\endaligned \leqno (2)$$ and
$$ \aligned
& \beta_n^{(h,r)}( q |\alpha_1 , \alpha_2, \cdots , \alpha_r) \\
& =\undersetbrace\text{$r$ times}\to{\int_{\Z}\cdots\int_{\Z}} [
\alpha_1x_1 +\cdots +\alpha_r x_r]^n q^{\sum_{i=1}^{r} (h-i)x_i}
d\mu_q(x_1)  \cdots d\mu_q(x_r). \endaligned \leqno (3) $$ These
can be written as
$$ \beta_n^{(h, r)}(w, q |\alpha_1 , \alpha_2, \cdots ,
 \alpha_r)= \sum_{j=0}^n  \binom nj q^{w_j} \beta_j^{(h, r)}( q |\alpha_1 , \alpha_2, \cdots ,
 \alpha_r)[w]^{n-j}. $$
By (3), we have
$$ \aligned
& \beta_n^{(h,r)}( q |\alpha_1 , \alpha_2, \cdots , \alpha_r) \\
& =\undersetbrace\text{$r$ times}\to{\int_{\Z}\cdots\int_{\Z}} [
\alpha_1x_1 +\cdots +\alpha_r x_r]^n q^{\sum_{i=1}^{r} (h-i)x_i}
d\mu_q(x_1)  \cdots d\mu_q(x_r)\\ & = \frac{1}{(1-q)^n}
\sum_{j=0}^n \binom nj (-1)^j \frac{(j \alpha_1 +h)(j \alpha_2
+h-1) \cdots (j \alpha_r +h-r+1)} {[j \alpha_1 +h][ j\alpha_2
+h-1] \cdots [j \alpha_r +h-r+1] }.
 \endaligned \leqno (4)$$
 Therefore we obtain the following :

 \proclaim{Theorem  1 }
 For any  $n \geq 0,$ we have
 $$\beta_n^{(h, r)}(w, q |\alpha_1 , \alpha_2, \cdots ,
 \alpha_r)  = \frac{1}{(1-q)^n}
\sum_{j=0}^n \binom nj (-q^w)^j \prod_{l=1}^r
\left(\frac{j\alpha_l +h-l+1}{[j \alpha_l +h-l+1] }\right).$$
 \endproclaim

Remark 1.  Note that
$$ \aligned
& \beta_0^{(2,1)}( q |1)=\frac 2{[2]}, \quad  \beta_1^{(2,1)}( q |1)=\frac {2q+1}{[2][3]}, \quad
\beta_2^{(2,1)}( q |1)= \frac {2q^2}{[3][4]},\\
& \beta_3^{(2,1)}( q |1)= - \frac {q^2(q-1)(2[3]+q)}{[3][4][5]}, \cdots ,  \\
&\beta_0^{(h,1)}( q |1)=\frac h{[h]}, \\
& \beta_1^{(h,1)}( q |1)= -\frac {(1+q+ \cdots + q^{h-1})+q(1+q+
\cdots + q^{h-2})+
\cdots +q^{h-1} }{[h][h+1]} \cdots \\
&\beta_0^{(2,2)}( q |1,1)=\frac {2!}{[2][1]}, \quad
\beta_1^{(2,2)}( q |1,1)=-\frac {2(q+2)}{[2][3]}, \\
&\beta_2^{(2,2)}( q |1,1)=-\frac {2((q-1)^2+ 5q)}{[3][4]}, \cdots  , \\
&\beta_0^{(r,r)}( q |\undersetbrace\text{$r$ times}\to{1, \cdots ,
1})=-\frac {r!}{[r][r-1]\cdots [2][1]},  \text{  ( cf. [4, 5, 6,
7, 8, 9, 10,  11] ) }.
 \endaligned $$
Remark 2. By the definition of $\beta_n^{(h,r)}( q |\alpha_1 ,
\alpha_2, \cdots , \alpha_r) ,$ we note that
$$\beta_n^{(h, r)}(0, q |\alpha_1 , \alpha_2, \cdots ,
 \alpha_r)=\beta_n^{(h,r)}( q |\alpha_1 , \alpha_2, \cdots , \alpha_r) .$$

 By (2),(3), it is easy to  see that
$$ \aligned
& \undersetbrace\text{$r$ times}\to{\int_{\Z}\cdots\int_{\Z}} [
\alpha_1x_1 +\cdots +\alpha_r x_r]^n q^{\sum_{i=1}^{r} (h-i)x_i}
d\mu_q(x_1)  \cdots d\mu_q(x_r)\\
& = (q-1)\undersetbrace\text{$r$
times}\to{\int_{\Z}\cdots\int_{\Z}} [ \alpha_1x_1 +\cdots
+\alpha_r x_r]^{n+1} q^{\sum_{i=1}^{r} (h-\alpha_i-i)x_i}
d\mu_q(x_1) \cdots
d\mu_q(x_r) \\
& + \undersetbrace\text{$r$ times}\to{\int_{\Z}\cdots\int_{\Z}} [
\alpha_1x_1 +\cdots +\alpha_r x_r]^{n} q^{\sum_{i=1}^{r}
(h-\alpha_i-i)x_i} d\mu_q(x_1) \cdots d\mu_q(x_r) .
 \endaligned $$
Thus we have
$$ \beta_m^{(h, r)}( q |\undersetbrace\text{$r$
times}\to{ 1 , 1, \cdots , 1})= (q-1)\beta_m^{(h-1, r)}( q
|\undersetbrace\text{$r$ times}\to{ 1 , 1, \cdots , 1}) +
\beta_m^{(h-1, r)}( q |\undersetbrace\text{$r$ times}\to{ 1 , 1,
\cdots , 1}).$$ It is easy to see that
$$ \aligned
& \undersetbrace\text{$r$ times}\to{\int_{X}\cdots\int_{X}} [ w+
\alpha_1x_1 +\cdots +\alpha_r x_r]^n q^{\sum_{i=1}^{r} (h-i)x_i}
d\mu_q(x_1)  \cdots d\mu_q(x_r)\\
& = [d]^{n-r} \sum_{i_1,i_2, \cdots i_r=0}^{d-1} q^{(i_1+ \cdots
i_r)h-i_2-2i_3- \cdots -(r-1)i_r} \\ & \times
\undersetbrace\text{$r$ times}\to{\int_{\Z}\cdots\int_{\Z}} [
\frac{w+ \alpha_1i_1 +\cdots +\alpha_r i_r}{d} + \alpha_1x_1 +
\cdots + \alpha_r x_r; q^d ]^{n} \\ &  \times   q^{x_1(h-1)d+
\cdot + x_r(h-r)d } d\mu_{q^d}(x_1) \cdots d\mu_{q^d}(x_r) .
 \endaligned  \leqno (5) $$
 From (3), (5), we have the following:

\proclaim{Theorem  2 }
 For any  positive integer $n ,$ we have
 $$
\aligned
  \beta_n^{(h, r)}(w, q |\alpha_1 , \alpha_2, \cdots ,
 \alpha_r)
  &=[d]^{n-r} \sum_{i_1,i_2, \cdots i_r=0}^{d-1} q^{(i_1+ \cdots
i_r)h-i_2-2i_3- \cdots -(r-1)i_r} \\ & \times
   \beta_n^{(h,r)}(
\frac{w+\alpha_1 i_1 + \cdots + \alpha_r i_r }{d}, q^d | \alpha_1,
\cdots , \alpha_r).
 \endaligned $$
 Moreover,
$$ \aligned
 \beta_n^{(h,r)}(wd, q |\alpha_1 , \alpha_2, \cdots , \alpha_r)
& =[d]^{n-r} \sum_{i_1,i_2, \cdots i_r=0}^{d-1} q^{(i_1+ \cdots
i_r)h-i_2-2i_3- \cdots -(r-1)i_r} \\ & \times
   \beta_n^{(h,r)}(
w+ \frac{\alpha_1 i_1 + \cdots + \alpha_r i_r }{d}, q^d |
\alpha_1, \cdots , \alpha_r).
\endaligned $$
\endproclaim

Remark 3. Note that
$$ \lim_{q \rightarrow 1} \beta_n^{(h,r)}(w, q |\alpha_1 , \alpha_2, \cdots , \alpha_r)
=B_n^{(r)}(w|\alpha_1 , \alpha_2, \cdots , \alpha_r), \text{ ( see
Eq. (1) ) } .$$
 Hence,
$\beta_n^{(h,r)}(w, q |\alpha_1 , \alpha_2, \cdots , \alpha_r)$
can be considered by the $q$-analogue of Barnes' multiple
Bernoulli polynomials. Now, we will give the inverse formula of
Eq. (4).

Indeed we see
$$ \aligned
& \sum_{i=0}^n \binom ni (q-1)^i \beta_i^{(h,r)}( q |\alpha_1 , \alpha_2, \cdots , \alpha_r)     \\
&= \sum_{i=0}^n \binom ni (q-1)^i \undersetbrace\text{$r$
times}\to{\int_{\Z}\cdots\int_{\Z}} [ \alpha_1x_1 +\cdots
+\alpha_r x_r]^i q^{\sum_{j=1}^{r} (h-j)x_j}
d\mu_q(x_1)  \cdots d\mu_q(x_r) \\
& = \undersetbrace\text{$r$ times}\to{\int_{\Z}\cdots\int_{\Z}}
q^{n( \alpha_1x_1 +\cdots +\alpha_r x_r)} q^{\sum_{j=1}^{r}
(h-j)x_j} d\mu_q(x_1)  \cdots d\mu_q(x_r) \\
&= \frac{(n \alpha_1 + h)(n \alpha_2 + h-1) \cdots (n \alpha_r +
h-r+1)}{[n \alpha_1 + h][n \alpha_2 + h-1] \cdots [n \alpha_r +
h-r+1]}.
 \endaligned $$
 Therefore we obtain the following:
\proclaim{Theorem  3 }
 For  $h \in \Bbb Z_{\geq 0},$ we have
 $$
 \sum_{i=0}^m \binom mi (q-1)^i \beta_i^{(h,r)}( q |\alpha_1 , \alpha_2, \cdots , \alpha_r)
 = \prod_{j=1}^r\left( \frac{m \alpha_j +h-j+1}{[m \alpha_j +h-j+1]} \right).$$
 \endproclaim

 Let $\chi$ be a  Dirichlet character  with
conductor $d \in \Bbb Z_{\geq 0}$. Then we define the generalized
Changhee  $q$-Bernoulli numbers  as follows: For $m \geq 0$,
$$
\aligned
 & \beta_{m,\chi}^{(h, r)}(q |\alpha_1 , \alpha_2, \cdots ,
 \alpha_r) \\ & = \undersetbrace\text{$r$ times}\to{\int_{X}\cdots\int_{X}} [
\alpha_1x_1 +\cdots +\alpha_r x_r]^m q^{\sum_{j=1}^{r} (h-j)x_j}
\left( \prod_{j=1}^r \chi(x_j)\right)d\mu_q(x_1)  \cdots
d\mu_q(x_r)
 \endaligned  . \leqno (6) $$
By simple calculation, we see that
$$
\aligned
 &  \undersetbrace\text{$r$ times}\to{\int_{X}\cdots\int_{X}} [
\alpha_1x_1 +\cdots +\alpha_r x_r]^m q^{\sum_{j=1}^{r} (h-j)x_j}
\left(\prod_{j=1}^r \chi(x_j)\right)d\mu_q(x_1)  \cdots d\mu_q(x_r) \\
&=[d]^{m-r} \sum_{i_1, \cdots, i_r=0}^{d-1} q^{ hi_1+ \cdots +
(h-r+1)i_r} \left( \prod_{j=1}^r \chi(x_j)\right) \\
& \times \undersetbrace\text{$r$
times}\to{\int_{\Z}\cdots\int_{\Z}} [ \frac{ \alpha_1i_1 +\cdots
+\alpha_r i_r}{d} + \alpha_1x_1 + \cdots + \alpha_r x_r; q^d ]^{m}
\\ & \times q^{\sum_{j=1}^{r}x_j (h-j)d} d\mu_{q^d}(x_1)  \cdots
d\mu_{q^d}(x_r).
 \endaligned  \leqno (7)$$
 By (2),(7), we have the following :
\proclaim{Theorem  4 }
 For  $h \in \Bbb Z_{\geq 0},$ we have
 $$
\aligned
 \beta_{m,\chi}^{(h, r)}(q |\alpha_1 , \alpha_2, \cdots
, \alpha_r)    &= [d]^{m-r} \sum_{i_1, \cdots, i_r=0}^{d-1} q^{
hi_1+ \cdots +
(h-r+1)i_r} \\
&\times \left( \prod_{j=1}^r \chi(x_j)\right)\beta_{m}^{(h,
r)}(\frac{\alpha_1 i_1+ \cdots + \alpha_r i_r}{d},q^d |\alpha_1 ,
\alpha_2, \cdots , \alpha_r).
 \endaligned $$
 \endproclaim

Remark.  By using our formulae in the case of $q=1$, we  can
obtain many new formulas which are related to the multiple Barnes'
Bernoulli numbers, cf.[1, 2, 4,  8,  10, 11].

ACKNOWLEDGEMENTS: This paper was supported by Korea Research
Foundation Grant((KRF-2003-05-C00009) .

\enddocument